\documentclass[10pt, conference]{IEEEtran}
\usepackage[english]{babel }

\usepackage[margin=54pt]{geometry}

\title{\vspace*{18pt} Distributed Primary Frequency Control through Multi-Terminal HVDC Transmission Systems}

\author{ \IEEEauthorblockA{Martin Andreasson$^{\IEEEauthorrefmark{1}}$, Roger Wiget, Dimos V. Dimarogonas, Karl H. Johansson and G\"oran Andersson
}
\\
\thanks{Martin Andreasson, Dimos V. Dimarogonas and Karl H. Johansson are with the ACCESS Linnaeus Centre, KTH Royal Institute of Technology, Stockholm, Sweden. Roger Wiget and G\"oran Andersson are with the Power Systems Laboratory, ETH Zurich, Zurich, Switzerland. This work was supported in part by the European Commission by the Hycon2 project, the Swedish Research Council (VR) and the Knut and Alice Wallenberg Foundation. 
We would like to thank the anonymous reviewers for their valuable comments.
\IEEEauthorrefmark{1} Corresponding author. E-mail: mandreas@kth.se}
}
\IEEEoverridecommandlockouts

\usepackage{mathrsfs}
\usepackage{amsfonts}
\usepackage{amssymb}
\usepackage{amsmath}
\usepackage{graphicx}
\usepackage{amsthm}
\usepackage{commath}
\usepackage{tikz}
\usepackage{tkz-graph}
\usepackage[europeanresistors]{circuitikz}
\usetikzlibrary{arrows}
\usepackage{flushend}
\usepackage{psfrag}
\usepackage{pgfplots}
\usepackage{todonotes}
\usepackage{blockdiagram}
\usetikzlibrary{external}
\usepackage{booktabs}
\newcommand{\raa}[1]{\renewcommand{\arraystretch}{#1}}

\newtheorem{theorem}{Theorem}

\newtheorem{remark}{Remark}

\newtheorem{assumption}{Assumption}
\newtheorem{objective}{Objective}

\DeclareMathOperator*{\diag}{diag}

\newcommand{\beq}{\begin{equation}}
\newcommand{\eeq}{\end{equation}}
\newcommand{\bq}{\begin{eqnarray}}
\newcommand{\eq}{\end{eqnarray}}
\newcommand{\bqn}{\begin{eqnarray*}}
\newcommand{\eqn}{\end{eqnarray*}}
\newcommand{\bee}{\begin{enumerate}}
\newcommand{\eee}{\end{enumerate}}

\newlength\fheight
\newlength\fwidth

\begin{document}
\maketitle

\begin{abstract}
This paper presents a decentralized controller for sharing primary AC frequency control reserves through a multi-terminal HVDC grid.
By using Lyapunov arguments, the proposed controller is shown to stabilize the equilibrium of the closed-loop system consisting of the interconnected AC and HVDC grids, given any positive controller gains. The static control errors resulting from the proportional controller are quantified and bounded by analyzing the equilibrium of the closed-loop system. The proposed controller is applied to a test grid consisting of three asynchronous AC areas interconnected by an HVDC grid, and its effectiveness is validated through simulation.
\end{abstract}

\section{Introduction}

Transmitting power over long distances with minimal losses is one of the greatest challenges in today's power transmission systems. The strong rising share of renewables increased the distances between power generation and consumption. This is a driving factor behind long-distance power transmission. One such example are large-scale off-shore wind farms, which often require power to be transmitted in cables over long distances to the mainland power grid \cite{breseti2007HVDC}. High-voltage direct current (HVDC) power transmission is a commonly used technology for long-distance power transmission. Its higher investment costs compared to AC transmission lines are compensated by its lower resistive losses for sufficiently long distances \cite{melhem2013electricity}. The break-even point, i.e., the point where the total construction and operation costs of overhead HVDC and AC lines are equal, is typically 500-800 km \cite{padiyar1990hvdc}. However, for cables, the break-even point is typically less than 50 km \cite{Hertem2010technical}. Increased use of HVDC for electrical power transmission suggests that future HVDC transmission systems are likely to consist of multiple terminals connected by several HVDC transmission lines \cite{Haileselassie2013Power}. Such systems are referred to as Multi-terminal HVDC (MTDC) systems in the literature. The main technical obstacle to overcome in order to realize such MTDC is the development of a DC breaker \cite{Franck2011HVDC}. There are a few advanced ideas to realize this device in the near future \cite{callavik2012hybrid}.

Maintaining an adequate DC voltage is the single most important practical control problem for HVDC transmission systems. If the DC voltage deviates too far from the nominal operational voltage, equipment could be damaged, resulting in loss of power transmission capability and high costs.

Many existing AC transmission grids are connected through HVDC links, usually used for bulk power transfer between the AC areas. The fast operation of the DC converters however would also enable frequency regulation of one of the connected AC grids through the HVDC link. One practical example of this is the island of Gotland in Sweden, which is only connected to the main Nordic grid through an HVDC cable \cite{axelsson2001gotland}. However, since the main Nordic AC grid has orders of magnitudes higher inertia than the AC grid of Gotland, the influence of the frequency regulation on the main grid will be negligible. 

By connecting several AC grids by an MTDC system, primary frequency regulation reserves may be shared, which reduces the need for frequency regulation reserves in the individual AC systems \cite{li2008frequency}. In \cite{dai2010impact}, distributed control algorithms have been applied to share primary frequency control reserves of asynchronous AC transmission systems connected through an MTDC system. However, the proposed controller requires a slack bus to control the DC voltage, defeating the purpose of distributing the primary frequency regulation reserves. In \cite{Andreasson2014_IFAC, andreasson2014control}, distributed controllers for secondary voltage control of MTDC systems are proposed, which do not rely on a slack bus. Both of the aforementioned controllers however rely on the presence of a communication network. While a communication network might already be present, it introduces the issue of time delays, due to large geographical distances in MTDC systems, and has a certain outage risk. The impacts of delays have been analyzed in \cite{dai2010impact}, and have been found to seriously degrade performance and destabilize the power system.
A distributed controller without the need of a slack bus is proposed in \cite{dai2013voltage}. Stability of the equilibrium is guaranteed in the absence of communication delays. However, the voltage dynamics of the HVDC system are neglected. Moreover, the implementation of the controller is not realistic, as every local controller needs to access the DC voltages of all terminals.
In \cite{dai2011voltage} and \cite{silva2012provision}, decentralized controllers are employed to share primary frequency control reserves. In \cite{silva2012provision} no stability analysis of the closed-loop system is performed, whereas \cite{dai2011voltage} guarantees stability of the equilibrium provided that the connected AC areas have identical parameters. In \cite{taylordecentralized}, optimal decentralized controllers for AC systems connected by HVDC systems are derived. In all aforementioned references the voltage dynamics of the HVDC system are neglected.

Due to the inherent difficulties of time-delays, we propose a decentralized proportional controller for distributing primary frequency control reserves, which relies only on local measurements. The controller is shown to distribute the primary frequency control reserves between the connected AC systems, while maintaining an adequate DC voltage. In contrast to \cite{dai2011voltage}, we prove that the equilibrium of the closed-loop system is globally asymptotically stable for any set of system parameters and controller gains by using Lyapunov arguments. We also explicitly model the voltage dynamics of the MTDC system, and extend our result to AC systems consisting of multiple generators in simulations. Due to inherent properties of proportional controllers, the steady-state values of the voltages and frequencies will deviate from their reference values. We quantify these deviations by provable upper bounds.

 The remainder of this paper is organized as follows. In Section \ref{sec:prel}, the mathematical notation is defined. In Section \ref{sec:model}, the system model and the control objectives are defined. In Section \ref{sec:dec_control}, a decentralized proportional controller for distributing primary frequency control is analyzed. In Section \ref{sec:equilibrium}, the equilibrium of the closed-loop system is analyzed. 
 In Section \ref{sec:simulations}, simulations of the distributed controller on a four-terminal MTDC test system are provided, showing the effectiveness of the proposed controller. The paper ends with a discussion and concluding remarks in Section \ref{sec:discussion}.

\section{Preliminaries}
\label{sec:prel}
Let $\mathcal{G}$ be a graph. Denote by $\mathcal{V}=\{ 1,\hdots, n \}$ the vertex set of $\mathcal{G}$, and by $\mathcal{E}=\{ 1,\hdots, m \}$ the edge set of $\mathcal{G}$. Let $\mathcal{N}_i$ be the set of neighboring vertices to $i \in \mathcal{V}$.
Denote by $\mathcal{B}$ the vertex-edge incidence matrix of $\mathcal{G}$, and let $\mathcal{\mathcal{L}_W}=\mathcal{B}W\mathcal{B}^T$ be the weighted Laplacian matrix of $\mathcal{G}$, with edge-weights given by the  elements of the diagonal matrix $W$. We denote the space of real-valued $n\times m$-valued matrices by $\mathbb{R}^{n\times m}$.
Let $\mathbb{C}^-$ denote the open left half complex plane, and $\bar{\mathbb{C}}^-$ its closure. We denote by $c_{n\times m}$ a vector or matrix of dimension $n\times m$, whose elements are all equal to $c$. For a symmetric matrix $A$, $A>0 \;(A\ge 0)$ is used to denote that $A$ is positive (semi) definite. $I_{n}$ denotes the identity matrix of dimension $n$. For simplicity, we will often drop the notion of time dependence of variables, i.e., $x(t)$ will be denoted $x$ for simplicity. Let $\norm{\cdot}_\infty$ denote the maximal absolute value of the elements of a vector.

\section{Model and problem setup}
\label{sec:model}
We will here give a unified model for an MTDC system interconnected with several asynchronous AC systems.
We consider an MTDC transmission system consisting of $n$ converters, each connecting to an AC system, denoted $1, \dots, n$. The converters are assumed to be connected by an MTDC transmission grid. The dynamics of converter $i$ is assumed to be given by
\begin{align}
\begin{aligned}
C_i \dot{V}_i &= -\sum_{j\in \mathcal{N}_i} \frac{1}{R_{ij}}(V_i -V_j) + I_i^{\text{inj}} ,
\end{aligned}
\label{eq:voltage}
\end{align}
where $V_i$ is the voltage of converter $i$, $C_i>0$ is its capacitance, $I_i^{\text{inj}}$ is the injected current from an AC grid connected to the DC converter.  The constant $R_{ij}$ denotes the resistance of the HVDC transmission line connecting the converters $i$ and $j$. 
The graph corresponding to the HVDC line connections is assumed to be connected. 
 The AC system is assumed to consist of a single generator which is connected to the corresponding DC converter, representing an aggregate model of the AC grid. The dynamics of the AC system are given by the swing equation \cite{machowski2008power}:
\begin{align}
m_i \dot{\omega}_i &=  -K_i^{\text{droop}} (\omega_i-\omega^{\text{ref}}) + P_i^\text{nom} + P_i^{{m}} - P_i^{\text{inj}}, \label{eq:frequency}
\end{align}
where $\omega_i$ is the frequency of the generator, $\omega^{\text{ref}}$ is the reference frequency and $m_i>0$ is its moment of inertia. The constant $P_i^\text{nom}$ is the nominal generated power of generator $i$, $P^m_i$ is the uncontrolled deviation from the nominal generated power, $P_i^{\text{inj}}$ is the power injected to the DC system through the converter and $K_i^{\text{droop}}>0$ is the gain of the frequency droop controller of the generator.
We define $P^\text{droop}_i=-K_i^{\text{droop}} (\omega_i-\omega^{\text{ref}})$, and state the control objective.
\begin{objective}
\label{obj:1}
The primary frequency control action should be distributed fairly amongst the generators, i.e.
\begin{align*}
\lim_{t\rightarrow \infty } \left|  P_i^{\text{droop}}(t) + \frac{1}{n} \sum_{i=1}^n P_i^m  \right| \le  e^{\text{droop}} \quad \forall i = 1, \dots, n,
\end{align*}
where $e^{\text{droop}}$ is a given scalar.
Furthermore, the frequencies of the AC systems, as well as the converter voltages, should not deviate too far from their nominal values, i.e.
\begin{align*}
\lim_{t\rightarrow \infty} |V_i(t)-V_i^{\text{ref}}| &\le e^{{V} } \quad \forall i = 1, \dots, n \\
\lim_{t\rightarrow \infty} |\omega_i(t)-\omega^{\text{ref}}| &\le e^{{\omega} } \quad \forall i = 1, \dots, n, \\
\end{align*}
where $V_i^{\text{ref}}$ is the reference DC voltage of converter $i$, $\omega^{\text{ref}}$ is the reference frequency and $e^{{V}}$ and $e^{{\omega} }$ are given scalars.
\end{objective}

\section{Decentralized MTDC control}
\label{sec:dec_control}
In this section we propose a decentralized controller for the frequency control of AC systems connected through an MTDC network. This controller does not rely on a single voltage regulator for the MTDC system, but the voltage regulation is distributed among all converters.
The local controller governing the power injections into the MTDC network is given by
\begin{align}
\label{eq:voltage_control}
\begin{aligned}
P_i^{\text{inj}} = P_i^{\text{inj, nom}} + K_i^{{\omega}} (\omega_i - \omega^{\text{ref}}) + K_i^{{V}}(V_i^{\text{ref}}-V_i),
\end{aligned}
\end{align}
where $P_i^{\text{inj, nom}}$ is the nominal injected power, and $K_i^{{\omega}}>0$ and $K_i^{{V}}>0$ are positive controller gains for all $i=1, \dots, n$. The HVDC converter is assumed to be perfect and instantaneous, i.e., injected power on the AC side is immediately converted to DC power without losses. Furthermore the dynamics of the converter are ignored, implying that the converter tracks the output of controller \eqref{eq:voltage_control} perfectly. This assumption is reasonable due to the dynamics of the converter typically being orders of magnitudes faster than the AC dynamics.
The relation between the injected HVDC current and the injected AC power is thus given by
\begin{align}
V_iI_i^{\text{inj}} = P_i^{\text{inj}}. \label{eq:power-current_nonlinear}
\end{align}
By assuming that all voltages are at the same nominal value, i.e., $V_i=V^{\text{nom}}$ for all $i=1, \dots, n$ in the above equation, the following linear relation is obtained
\begin{align} 
V^{\text{nom}}I_i^{\text{inj}} = P_i^{\text{inj}}. \label{eq:power-current}
\end{align}
Combining the voltage dynamics \eqref{eq:voltage}, the frequency dynamics \eqref{eq:frequency}, the voltage controller \eqref{eq:voltage_control} and the power-current relationship \eqref{eq:power-current}, we obtain the following closed-loop dynamics
\begin{align}
\begin{bmatrix}
\dot{\omega} \\ \dot{V}
\end{bmatrix}
&= \underbrace{\begin{bmatrix}
-M(K^\omega + K^{\text{droop}}) & MK^V \\
\frac{1}{V^{\text{nom}}}EK^\omega & -E\left(\mathcal{L}_R + \frac{K^V}{V^{\text{nom}}} \right)
\end{bmatrix}}_{\triangleq A}
\begin{bmatrix}
\omega \\ V
\end{bmatrix} \nonumber \\
&+ \begin{bmatrix}
M\left((K^\omega + K^{\text{droop}}) \omega^{\text{ref}}1_{n\times 1} - K^V V^{\text{ref}} \right) \\
E\left(\frac{1}{V^\text{nom}} K^V V^\text{ref} -\frac{\omega^{\text{ref}}}{{V^{\text{nom}}}} K^\omega1_{n\times 1} \right)
\end{bmatrix}  \nonumber \\
&+
\begin{bmatrix}
M(P^m + P^\text{nom}-P^{\text{inj, nom}}) \\
\frac{1}{V^{\text{nom}}} E  P^{\text{inj, nom}}
\end{bmatrix}
\label{eq:cl_dynamics_vec}
\end{align}
where
$ \omega = [\omega_1, \dots, \omega_n]^T$, 
$ V =[V_1, \dots, V_n]^T$, 
 $M=\diag({m_1}^{-1}, \hdots , {m_n}^{-1})$ is a matrix of inverse generator inertia, 
 $E=\diag([C_1^{-1}, \dots, C_n^{-1}])$ is a matrix of electrical elastances, 
 $K^\omega = \diag([K^\omega_1, \dots, K^\omega_n])$, 
 $K^{\text{droop}} = \diag([K^{\text{droop}}_1, \dots$, $K^{\text{droop}}_n])$,
  $K^V = \diag([K^V_1, \dots, K^V_n])$, 
  $P^\text{nom} =[P^\text{nom}_1,\hdots, P^\text{nom}_n]^T$, 
  $P^\text{inj, nom} =[P^\text{inj, nom}_1,\hdots, P^\text{inj, nom}_n]^T$, $P^m =[P^m_1,\hdots, P^m_n]^T$,
 and $\mathcal{L}_R$ is the weighted Laplacian matrix of the graph representing the HVDC transmission lines, denoted $\mathcal{G}_R$, whose edge-weights are given by the conductances $\frac{1}{R_{ij}}$. The following assumption is made on the nominal generated power and the nominal injected power.

\begin{assumption}
\label{ass:balances_power}
$P^\text{nom}=P^\text{inj, nom}$. 
\end{assumption}
\begin{remark}
Assumption~\ref{ass:balances_power} implies that the reference frequency and reference voltages define an equilibrium of the closed-loop system when the deviation from the nominal power generation is zero.
\end{remark}
 We define the incremental frequencies and voltages as
\begin{align}
\hat{ \omega} &= \omega-\omega^{\text{ref}}1_{n\times 1} \label{eq:delta_omega}
 \\
\hat{ V} &= V- V^{ref} \label{eq:delta_V}.
\end{align}
By Assumption \ref{ass:balances_power}, the decentralized MTDC control system given by \eqref{eq:cl_dynamics_vec}, can be written as
\begin{align}
\begin{bmatrix}
\dot{\hat{\omega}} \\ \dot{\hat{V}}
\end{bmatrix}
&= {\begin{bmatrix}
-M(K^\omega + K^{\text{droop}}) & MK^V \\
\frac{1}{V^{\text{nom}}}EK^\omega & -E\left(\mathcal{L}_R + \frac{K^V}{V^{\text{nom}}} \right)
\end{bmatrix}}
\begin{bmatrix}
\hat{\omega} \\ \hat{V}
\end{bmatrix} \nonumber \\
&+
\begin{bmatrix}
M P^m \\
0_{n\times 1}
\end{bmatrix}. \label{eq:cl_dynamics_vec_delta}
\end{align}
Assume that the system matrix of \eqref{eq:cl_dynamics_vec_delta}, $A$, is full-rank, which ensures that unique equilibrium of \eqref{eq:cl_dynamics_vec_delta}  exists. Denote this equilibria $x_{0}=[\omega_{0}^T, V_{0}^T]^T$. Define $\bar{x}\triangleq [\bar{\omega}^T, \bar{V}^T]^T =[\hat{\omega}^T, \hat{V}^T]^T - [\omega_{0}^T, V_{0}^T]^T$. 
Now: 
\begin{align}
\dot{\bar{x}} = A \bar{x} \label{eq:dynamics_A_decentralized_shifted} 
\end{align}
with the origin as the unique equilibrium of the above dynamical system. We are now ready to show the main stability result of this section.

\begin{theorem}
\label{th:stability_passivity_1}
The equilibrium of the decentralized MTDC control system given by \eqref{eq:cl_dynamics_vec_delta} is globally asymptotically stable.
\end{theorem}
\begin{proof}
First consider the Lyapunov function candidate 
\begin{align}
W(\bar{\omega}, \bar{V}) &= \frac 12 \bar{\omega}^T K^\omega (K^V)^{-1} M^{-1}\bar{\omega} + \frac{V^\text{nom}}{2} \bar{V}^T C \bar{V}, \label{eq:lyap_hvdc_decentralized_projected}
\end{align}
where  $C=\diag([C_1, \dots, C_n])$. 
Clearly $W(\bar{\omega}, \bar{V})$ is positive definite and radially unbounded. Differentiating \eqref{eq:lyap_hvdc_decentralized_projected} with respect to time along trajectories of \eqref{eq:dynamics_A_decentralized_shifted}, we obtain
\begin{align*}
&\dot{W}(\bar{\omega}, \bar{V})  \\
&= \bar{\omega}^T K^\omega (K^V)^{-1} M^{-1}\dot{\bar{\omega}} + V^\text{nom} \bar{V}^T E \dot{\bar{V}} + \bar{\eta}' \dot{\bar{\eta}}' \\
&= \bar{\omega}^T \big( -K^\omega (K^V)^{-1}(K^\omega + K^\text{droop})\bar{\omega} + K^\omega \bar{V} \big) \\
&\;\;\;\; + \bar{V}^T \Big( K^\omega \bar{\omega} - (V^\text{nom}\mathcal{L}_R {+} K^V)\bar{V} \Big) \\
&= -\bar{\omega}^T \big( -K^\omega (K^V)^{-1}(K^\omega + K^\text{droop})\bar{\omega} \\
&\;\;\;\; +  2 \bar{\omega}^T K^\omega \bar{V} - \bar{V}^T (V^\text{nom}\mathcal{L}_R + K^V)\bar{V}   \\
&= - \begin{bmatrix}
\bar{\omega}^T & \bar{V}^T 
\end{bmatrix}
\underbrace{\begin{bmatrix}
K^\omega (K^V)^{-1}(K^\omega + K^\text{droop}) & -K^\omega \\
-K^\omega & K^V 
\end{bmatrix}}_{\triangleq Q_1}
\begin{bmatrix}
\bar{\omega} \\ \bar{V}
\end{bmatrix}.
\end{align*} 
Clearly $\dot{W}(\bar{\omega}, \bar{V})< 0$ iff the symmetric matrix $Q_1$ is positive definite. By applying the Schur complement condition for positive definiteness, $Q_1$ is positive definite iff
\begin{eqnarray*}
K^\omega (K^V)^{-1}(K^\omega + K^\text{droop}) - K^\omega (K_V)^{-1} K^\omega \\
= K^\omega (K^V)^{-1} K^\text{droop} > 0.
\end{eqnarray*}
Hence $Q_1$ is always positive definite, and thus $\dot{W}(\bar{\omega}, \bar{V}< 0$, which concludes the proof.
\end{proof}

\begin{remark}
Note that the equilibrium of  \eqref{eq:cl_dynamics_vec_delta}  being globally asymptotically stable implies that all the eigenvalues of $A$ are stable, which ensures that the previous assumption that $A$ is full rank, is valid. 
\end{remark}

\begin{remark}
Note that  Theorem~\ref{th:stability_passivity_1} only guarantees the stability of the equilibrium. It does however not guarantee that Objective~\ref{obj:1} is fulfilled. 
\end{remark}

\section{Equilibrium analysis}
\label{sec:equilibrium}
We will now study the globally asymptotically stable equilibrium of \eqref{eq:cl_dynamics_vec_delta}, in order to bound the asymptotic voltage and frequency deviations from the reference values. We will furthermore show that the generated power in the AC grids will be shared fairly amongst the generators. 
 We make the following additional assumptions on the controller gains, in order to draw conclusions about the equilibrium of \eqref{eq:cl_dynamics_vec_delta}.
\begin{assumption}
\label{ass:scalar_1} The controller gains satisfy
$K^\omega_i=k^\omega, K^\text{droop}_i=k^\text{droop}, K^V_i=k^V \; \forall i=1, \dots, n$.
\end{assumption}

With the previous assumptions made, and having provided necessary stability conditions for the closed-loop system \eqref{eq:cl_dynamics_vec}, we are ready to analyze its equilibrium.
\begin{theorem}
\label{th:equilibrium}
Assume that Assumptions~\ref{ass:balances_power} and \ref{ass:scalar_1} hold, then Objective \ref{obj:1} is satisfied for the following coefficients
\begin{align*}
e^{\text{gen}} &=\frac{k^\text{droop}\max_i P^m_i}{k^\text{droop}+k^\omega} \left( (n-1) + \frac{k^V}{V^\text{nom}} \sum_{i=2}^n \frac{1}{\lambda_i(\mathcal{L}_R)} \right) \\
e^V &= \frac{k^\omega \left|1_{1\times n}P^m\right|}{nk^\text{droop}k^V}  +  \frac{k^\omega \max_i \left| P^m_i  \right| }{(k^\omega + k^\text{droop})V^\text{nom}} \sum_{i=2}^n \frac{1}{\lambda_i(\mathcal{L}_R)} \\
e^\omega &= \frac{1}{n k^\text{droop}} \left| \sum_{i=1}^n P^m_i \right| \\
&\;\;\;\; + \frac{\max_i |P^m_i|}{k^\text{droop}+k^\omega} \left( (n-1) + \frac{k^V}{V^\text{nom}} \sum_{i=2}^n \frac{1}{\lambda_i(\mathcal{L}_R)} \right).
\end{align*}
\end{theorem}

\begin{remark}
The error bounds $e^{\text{droop}}$ and $e^\omega$ can simultaneously be made arbitrarily small by choosing appropriate controller gains. However, the voltage error bound $e^V$ is lower bounded by a constant. This is of course due to the necessity of a relative voltage drop for having a power flow between in an HVDC line.
\end{remark}

\begin{proof}
Consider the equilibrium of \eqref{eq:cl_dynamics_vec}. Let  $\hat{\omega}$ and $\hat{V}$ be defined by \eqref{eq:delta_omega} -- \eqref{eq:delta_V}. By Assumption~\ref{ass:balances_power}, we obtain the following expression
\begin{align}
\label{eq:eq_delta_coordinates}
\begin{bmatrix}
-(K^\omega+K^\text{droop}) & K^V \\
K^\omega & -(K^V+ V^\text{nom} \mathcal{L}_R)
\end{bmatrix}
\begin{bmatrix}
\hat{\omega}  \\
\hat{V}
\end{bmatrix}
&=
\begin{bmatrix}
-P^m \\
0_{n\times 1}
\end{bmatrix}.
\end{align}
By multiplying the last $n$ rows of \eqref{eq:eq_delta_coordinates} with $\frac{k^\omega+ k^\text{droop}}{k^\omega}$ and adding to the first $n$ rows of \eqref{eq:eq_delta_coordinates}, we obtain by Assumption~\ref{ass:scalar_1}
\begin{align}
\label{eq:Delta_V_eq}
\underbrace{\left( \frac{(k^\omega+k^\text{droop})V^\text{nom}}{k^\omega} \mathcal{L}_R + \frac{k^\text{droop}k^V}{k^\omega}I_n \right)}_{\triangleq A_1}\hat{V} = P^m.
\end{align}
We write $\hat{V} = \sum_{i=1}^n a^1_i v^1_i$, where $v^1_i$ is the $i$th eigenvector of $A_1$, with the corresponding eigenvalue $\lambda^1_i$. Note that the coefficients $a^1_i$ are unique, since $A_1$ is symmetric, implying that its eigenvectors form an orthonormal basis. Substituting the eigenvector decomposition of $\hat{V}$ in \eqref{eq:Delta_V_eq} yields
\begin{align*}
A_1\hat{V} = A_1 \sum_{i=1}^n a^1_i v^1_i = \sum_{i=1}^n \lambda^1_i a^1_i v^1_i = P^m,
\end{align*}
which implies
\begin{align*}
a^1_i=\frac{(v^1_i)^TP^m}{\lambda^1_i}.
\end{align*}
Let the eigenvalues be ordered by their increasing values. Clearly $\lambda^1_1=  \frac{k^\text{droop}k^V}{k^\omega}$ and $v^1_1=\frac{1}{\sqrt{n}} 1_{n\times 1}$. This implies
\begin{align}
\label{eq:Delta_V_ss}
\hat{V} = \frac{k^\omega 1_{1\times n}P^m}{nk^\text{droop}k^V} 1_{n\times 1} + \sum_{i=2}^n \frac{(v^1_i)^TP^m}{\lambda^1_i}v^1_i.
\end{align}
By noting that
\begin{align}
\begin{aligned}
\lambda^1_i &= \frac{(k^\omega + k^\text{droop})V^\text{nom} \lambda_i(\mathcal{L}_R)+k^\text{droop}k^V}{k^\omega} \\
&\ge \frac{(k^\omega + k^\text{droop})V^\text{nom} \lambda_i(\mathcal{L}_R)}{k^\omega},
\end{aligned}
\label{eq:lambda_1_bound}
\end{align}
where $\lambda_i(\mathcal{L}_R)$ is the $i$th eigenvalue of $\mathcal{L}_R$,
we obtain the following bound on $\hat{V}$
\begin{align*}
\norm{\hat{V}}_\infty &\le  \frac{k^\omega \left|1_{1\times n}P^m\right|}{nk^\text{droop}k^V}  +  \frac{\max_i \left| P^m_i  \right| k^\omega}{(k^\omega + k^\text{droop})V^\text{nom}} \sum_{i=2}^n \frac{1}{\lambda_i(\mathcal{L}_R)} \\
& \le  \frac{k^\omega \left|\sum_{i=1}^nP^m_i\right|}{nk^\text{droop}k^V}  +  \frac{\max_i \left| P^m_i  \right|}{V^\text{nom}} \sum_{i=2}^n \frac{1}{\lambda_i(\mathcal{L}_R)}.
\end{align*}
From the first $n$ rows of \eqref{eq:eq_delta_coordinates}, and by substituting the expression for $\hat{V}$ from \eqref{eq:Delta_V_ss}, we have
\begin{align}
\begin{aligned}
\hat{\omega} &= \frac{k^V\hat{V} + P^m}{k^\omega+ k^\text{droop}} = \frac{1}{k^\omega + k^\text{droop}} \Bigg( \frac{k^\omega 1_{1\times n}P^m}{nk^\text{droop}} 1_{n\times 1}   \\
  &\;\;\;\;+  \sum_{i=2}^n \frac{k^V(v^1_i)^TP^m}{\lambda^1_i}v^1_i + P^m \Bigg).
\end{aligned}
\label{eq:Delta_omega}
\end{align}
By using the bound on $\lambda^1_i$ from \eqref{eq:lambda_1_bound}, we obtain
\begin{align*}
\norm{\hat{\omega}}_\infty &\le \frac{1}{k^\text{droop}} \Bigg( \frac{\left|\sum_{i=1}^n P^m_i\right|}{n}  \\
& \;\;\;\; +  \max_i |P^m_i|\Bigg( 1 + \frac{k^V}{V^\text{nom}} \sum_{i=2}^n \frac{1}{\lambda_i(\mathcal{L}_R)} \Bigg) \Bigg).
\end{align*}
Consider now the power generated by the voltage droop controller. By \eqref{eq:Delta_omega} we obtain
\begin{align*}
&{P^\text{droop} + \frac{1}{n} \sum_{i=1}^n P_i^m 1_{n\times 1}} = {-k^\text{droop}\hat{\omega} + \frac{1_{1\times n} P^m}{n} 1_{n\times 1} } \\
&= \frac{k^\text{droop}}{k^\omega+k^\text{droop}} \Bigg(- \frac{1}{n}\sum_{i=1}^n {P^m_i}1_{n\times 1}  + P^m  \\
&\;\;\;\;\;{+}  \sum_{i=2}^n \frac{k^V(v^1_i)^TP^m}{\lambda^1_i}v^1_i \Bigg) .
\end{align*}
By using the bound on $\lambda^1_i$ in \eqref{eq:lambda_1_bound}, we obtain
\begin{align*}
&\norm{{P^\text{droop} + \frac{1}{n} \sum_{i=1}^n P_i^m 1_{n\times 1}}}_\infty \le  \frac{k^\text{droop}}{k^\omega+k^\text{droop}} \Bigg(  \\
&\;\;\;\; \max_i \left| P^m_i \right|\left( 1+ \sum_{i=2}^n \frac{k^\omega}{(k^\omega + k^\text{droop})V^\text{nom} \lambda_i(\mathcal{L}_R)} \right) \Bigg) \\
& \le \frac{k^\text{droop}}{k^\omega+k^\text{droop}}  \max_i \left| P^m_i \right|\Bigg( 1+ \frac{1}{V^\text{nom}}\sum_{i=2}^n \frac{1}{ \lambda_i(\mathcal{L}_R)}  \Bigg),
\end{align*}
which completes the proof.
\end{proof}

\section{Simulations}
\label{sec:simulations}
\begin{figure}[htb]
\centering
\includegraphics[width=1.0\columnwidth]{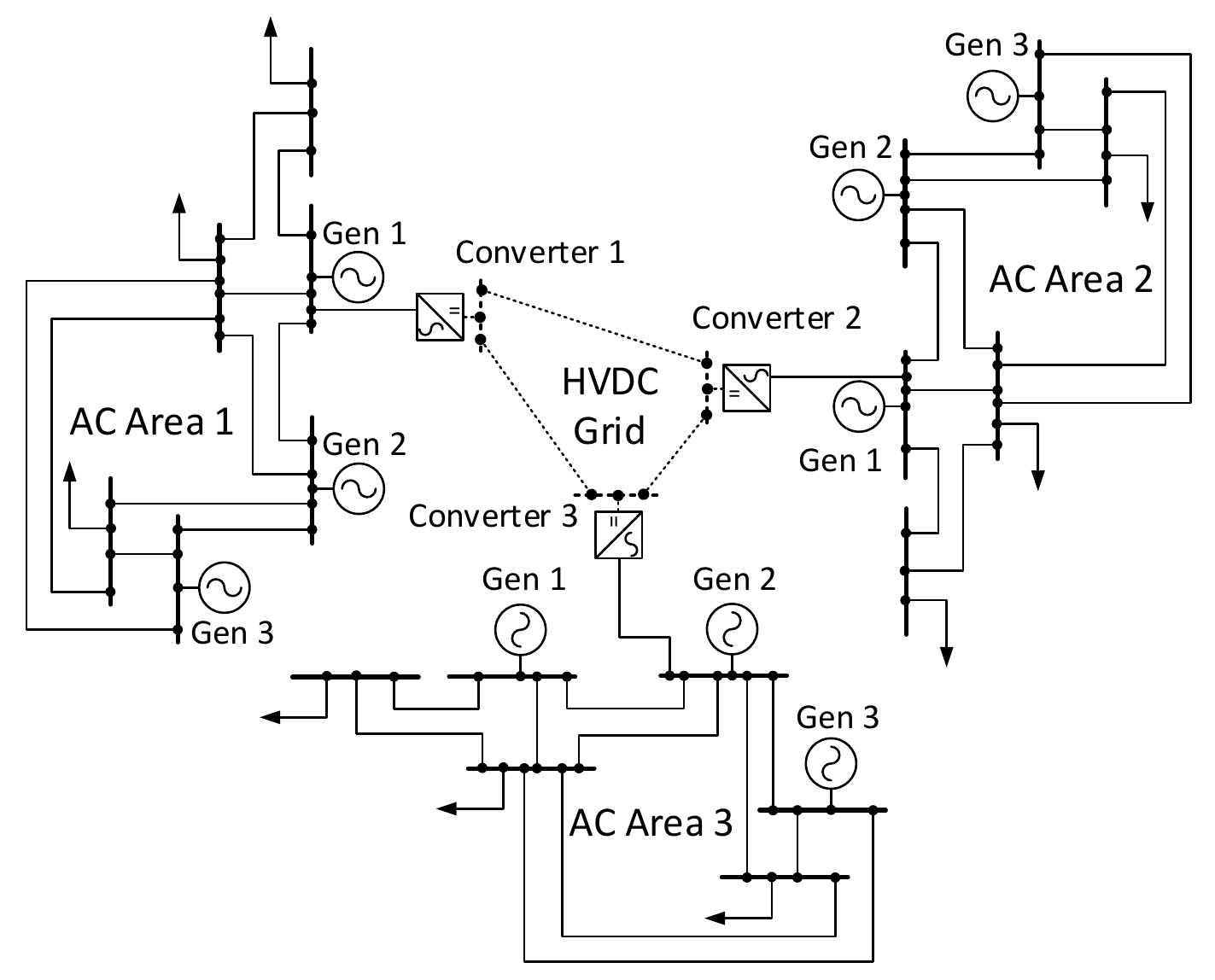}
\caption{Test grid consisting of 3 AC areas, connected by an MTDC grid consisting of 3 converter stations and 3 DC lines.}
\label{fig:testgrid}
\end{figure}

\newlength\figureheight
\newlength\figurewidth
\setlength\figureheight{4.4cm}
\setlength\figurewidth{6.6cm}
In this section, we simulate the proposed controller on an MTDC grid connecting three asynchronous AC areas, whose main purpose is bulk power transfer between the AC areas. The test grid consists of tree 6 bus AC grids, described in detail in \cite{wollenberg2006power}, connected with a 3 bus MTDC grid. In Figure \ref{fig:testgrid}, the topology of the interconnected MTDC-AC grid is shown.

\begin{figure}[th]
\input{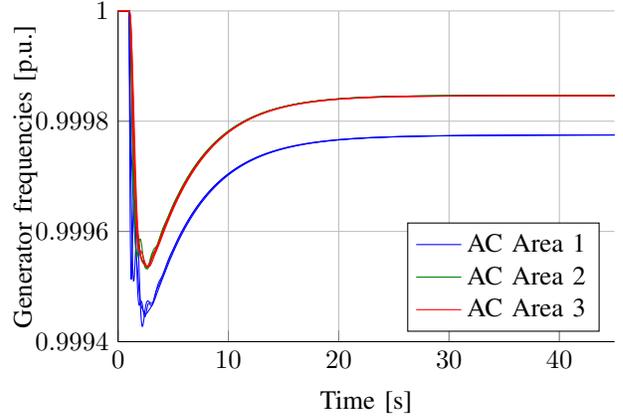}
\caption{Frequencies of the generator areas.}
\label{fig:generatorspeeds}
\end{figure}

\begin{figure}[th]
\center
\input{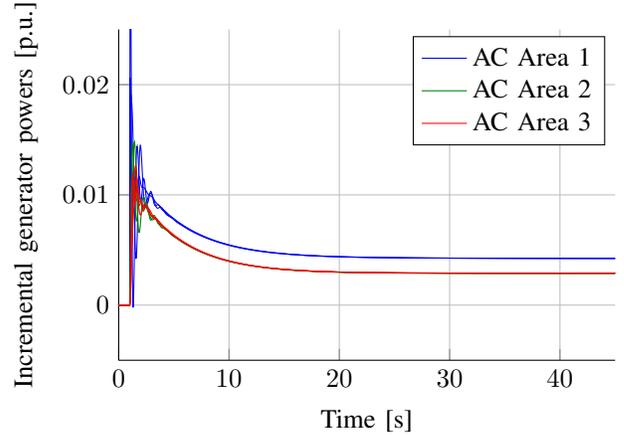}
\caption{Incremental generator power levels.}
\label{fig:GeneratorDelta}
\end{figure}

\begin{figure}[th]
\center
\input{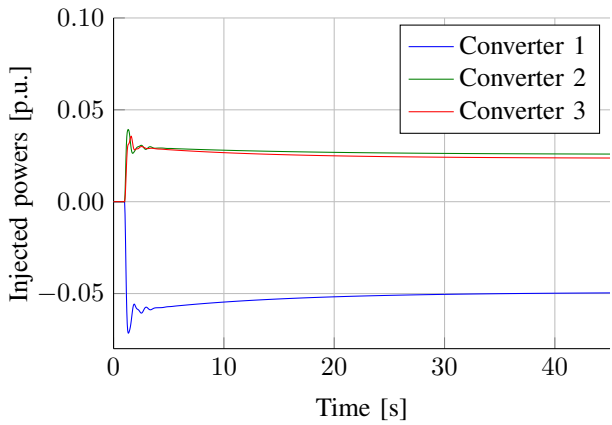}
\caption{Injected power levels at the converters.}
\label{fig:converterpower}
\end{figure}

\begin{figure}[th!]
\center
\input{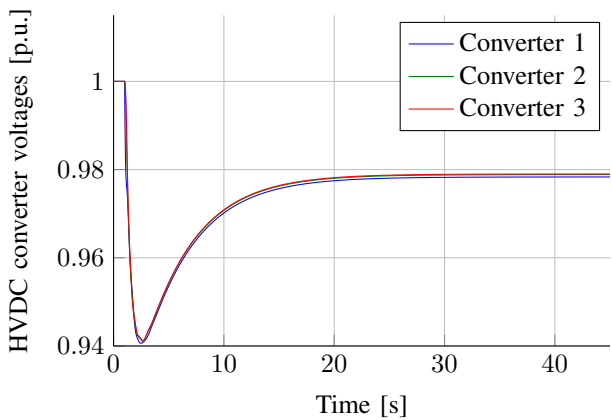}
\caption{Voltages of the DC converters.}
\label{fig:convertervoltages}
\end{figure}

\begin{table}
\centering
\caption{HVDC grid line parameters}
\label{tab:HVDCgridParameter}
\begin{tabular}{llll}\toprule
From & To & Resistance [p.u.]& Reactance [p.u.] \\ \midrule
  1 & 2 & 0.0015 & 0.01 \\
  1 & 3 & 0.0045 & 0.03 \\
  2 & 3 & 0.0015 & 0.01 \\
  \bottomrule
\end{tabular}
\end{table}
\begin{table}\centering
\raa{1.3}
\caption{Controller Parameter}
\label{tab:ControllerParameter}
\begin{tabular}{@{}llllll@{}}\toprule
$K^{\omega}_1$ &$K^{\omega}_2$&$K^{\omega}_3$& $K^\text{droop}_1$ &$K^\text{droop}_2$ &$K^\text{droop}_3$\\ \midrule
501&501& 501&667&667 &  667 \\
\bottomrule
\end{tabular}
\end{table}
Each converter station is controlled with \eqref{eq:voltage_control}.
While the converter dynamics are ignored due to their fast nature, the nonlinear relation \eqref{eq:power-current_nonlinear} is used to relate the injected AC powers with the HVDC currents.
 The physical system parameters and the controller parameters are given in Table \ref{tab:HVDCgridParameter}, \ref{tab:ControllerParameter}, respectively.
The simulation was conducted by using an extended version of MatDyn \cite{matsch}, taking also the HVDC dynamics into account. The simulation starts in steady-state, and at time $1$ s an immediate change in load from 0.7 p.u. (per-unit) to 0.8 p.u. occurs at bus 4 in the AC area 1. 
The local frequency controllers at the generators will react immediately to the resulting frequency drop, and start to accelerate. The frequencies of the generators are shown in Figure \ref{fig:generatorspeeds}.
After a few second all generator frequencies within the same AC area synchronize, and after about $30$ s the frequencies converge to the new equilibrium. The frequency deviation is larger in AC area $1$ than in the remaining AC areas, but the differences are rather small, in accordance with Theorem~\ref{th:equilibrium}.
Figure \ref{fig:GeneratorDelta} shows the changes in the power output of the generators. The disturbance is shared along all generators.
The injected powers through the converter are shown in Figure \ref{fig:converterpower}. Since the converter dynamics are much faster  than the AC systems, they are neglected in the simulation and it is assumed that the converter power tracks the controller output perfectly.
Due to the increased load, the DC voltages of all converters increase, see Figure \ref{fig:convertervoltages}. However, as predicted by Theorem~\ref{th:equilibrium}, both the absolute and relative voltage deviations are bounded.

\section{Discussion and Conclusions}
\label{sec:discussion}
In this paper we have proposed a decentralized proportional controller for sharing primary frequency control reserves in asynchronous AC systems connected through an MTDC system. The controller uses the local frequency in the AC grid and the local DC voltage as inputs in order to control the power injections into the MTDC grid. The resulting equilibria of the closed-loop system is shown to be globally asymptotically stable by using Lyapunov arguments, regardless of the controller parameters. It is also shown that the DC voltages and AC frequencies at the equilibrium are close to their nominal values. Furthermore the generated power from the primary frequency control is approximately shared fairly between the AC areas. The deviation from perfectly fair power sharing is quantified.
The proposed controller was simulated on a test system consisting of 3 AC areas combined with an MTDC grid to demonstrate its effectiveness. The paper constitutes a first step towards utilizing the increased flexibility which future MTDC grids will provide to the connected AC systems. Future work will focus on extending the primary proportional controller with secondary controllers, where communication and integral action will be necessary to eliminate static control errors. An extensive simulation study on more realistic grid topologies and dynamical models is also ongoing work. 
\bibliography{references}
\bibliographystyle{plain}
\end{document}